\documentclass[12pt]{article}
\usepackage{amsfonts,amssymb,amsmath,latexsym,amscd}
\newtheorem{theorem}{Theorem}[section]

\newtheorem{proposition}[theorem]{Proposition}
\newtheorem{corollary}[theorem]{Corollary}

\newcommand {\mat}      [1] {\left(\begin{array}{#1}}
\newcommand {\rix}          {\end{array}\right)}
\newcommand{\proof}{\par\noindent {\em Proof}. \ignorespaces}
\newcommand{\eproof}{\mbox{\hspace*{20pt} $\Box$}}
\newcommand{\BC}{\mathbb C}
\newcommand{\BR}{\mathbb R}

\newcommand{\eq} [1] {\begin{equation}\label{#1}}
\newcommand{\en} {\end{equation}}

\newcommand {\eqn}      {\begin{eqnarray}}
\newcommand {\enn}      {\end{eqnarray}}
\newcommand {\bstar}    {\begin{eqnarray*}}
\newcommand {\estar}    {\end{eqnarray*}}

\begin{document}

\title{On matrix polynomials with real roots}
\date{}

\author{Leonid Gurvits\thanks{Los Alamos National Laboratory, 
Los Alamos, NM 87545, USA. E-mail: gurvits@lanl.gov.}  \and
Leiba Rodman\thanks{College of William and Mary, Department
of Mathematics, P. O. Box 8795, Williamsburg, VA 23187-8795, USA. E-mail:
lxrodm@math.wm.edu. The research leading to this article 
was done while the second author visited LANL, whose hospitality is 
gratefully acknowledged. 
The research of second author is partially supported by NSF
Grant DMS-9988579.}}

\maketitle

\begin{abstract}
It is proved that 
the roots of combinations of matrix polynomials with real roots can be 
recast as
eigenvalues of combinations of real symmetric matrices, under certain 
hypotheses. The proof is based on recent solution of the
Lax conjecture. Several applications and corollaries, in particular
concerning hyperbolic matrix polynomials, are presented.  
\end{abstract}

\vspace{3mm}


\section{Main result}
\setcounter{equation}{0}

A polynomial is called {\em hyperbolic}
if all its roots are real. It is a classical well studied class of polynomials
(see, e.g., \cite{O}). There are at least two useful ways to 
extend this notion to polynomials with   
complex $n \times n$ matrix coefficients, in short, matrix polynomials. 
Thus, a monic (i.e., with leading coefficient $I_n$, the $n \times n$ identity 
matrix) matrix polynomial 
$ L(z)$
of degree $\ell$ is said to be {\em hyperbolic}, if for every nonzero 
$x\in \BC^n$, the $n$-dimensional 
vector space of columns with complex components, 
the polynomial equation 
\begin{equation}\label{12}  \langle L(z)x,x\rangle=0 \end{equation}
has $\ell$ real roots (counted with multiplicities). 
We denote here by $\langle \cdot, \cdot \rangle$ the standard inner product in 
$\BC^n$. An $n 
\times n$ monic matrix 
polynomial $L(z)$ 
of degree $\ell$ will be called {\em weakly hyperbolic} if 
${\rm det} L(z)=0$ has $n\ell$ real roots (multiplicities counted). 
Note that our terminology differs slightly from the terminology in some 
sources (for example, \cite{3}). 
Clearly, every hyperbolic matrix polynomial is weakly hyperbolic, and  
the coefficients of every hyperbolic matrix polynomial are 
Hermitian matrices. See, e.g., \cite{1,2,3,4} for the 
theory and applications of 
hyperbolic matrix and operator polynomials. 

In this note we prove the following theorem. It states that 
the roots of combinations of hyperbolic matrix polynomials  can be recast as
eigenvalues of combinations of real symmetric matrices, under certain 
hypotheses.
We denote by $\BR$ the field of real numbers. 

\begin{theorem}\label{th1} 
Let 
$$ L(z) =\sum_{j=0}^{\ell} L_jz^j, \qquad 
M(z)=\sum_{j=0}^{\ell} M_jz^j, \qquad M_{\ell}=L_{\ell}=I, $$
be two monic $n \times n$ matrix polynomials such that 
\begin{equation}\label{*} 
 \alpha L(z)+(1 -\alpha)M(z) \quad 
{\rm is} \ {\rm  weakly} \ {\rm  hyperbolic} \ {\rm  for} \ {\rm  every}\  \alpha 
\in  \BR. \end{equation} 
Assume in addition that the $n \times n$ matrix  
$L_{\ell-1}-M_{\ell-1} $
has $n$ real eigenvalues (counted with multiplicities). 
Then there exist $n\ell \times n\ell$ real symmetric matrices 
$A$ and $B$ such that 
for every $\alpha\in \BR$, the 
roots of ${\rm det}\, (\alpha L(z)+(1-\alpha)M(z))$, counted according to 
 their multiplicities, coincide with the eigenvalues of $
\alpha A+(1-\alpha)B$, also counted according to their multiplicities.

Conversely, if
the 
roots of ${\rm det}\, (\alpha L(z)+(1-\alpha)M(z))$, 
coincide with the eigenvalues of $
\alpha A+(1-\alpha)B$ $($counted with multiplicities$)$ for every 
$\alpha\in \BR$, where $A$ and $B$ are fixed real symmetric 
$n\ell \times n\ell$ matrices, then $(\ref{*})$ holds and  
all eigenvalues of $L_{\ell-1}-M_{\ell-1} $ are real.
\end{theorem}

\proof Let 
$$C_L=\mat{ccccc} 0 & I_n &0 & \ldots & 0 \\
0&0& I_n &\ldots &0 \\ \vdots & \vdots & \vdots & \ddots & \vdots \\
0&0&0& \ldots &I_n \\ -L_0& -L_1&-L_2& \ldots & -L_{\ell-1} \rix 
$$ and $$ C_M=
\mat{ccccc} 0 & I_n &0 & \ldots & 0 \\
0&0& I_n &\ldots &0 \\ \vdots & \vdots & \vdots & \ddots & \vdots \\
0&0&0& \ldots &I_n \\ -M_0& -M_1&-M_2& \ldots & -M_{\ell-1} \rix
$$
be the companion 
matrices of $L(z)$ and 
of $M(z)$, respectively. Then 
$\alpha C_L+(1-\alpha)C_M$ is the companion matrix of 
$\alpha L(z)+(1-\alpha)M(z)$, and therefore the roots of 
${\rm det}\, (\alpha L(z)+(1-\alpha)M(z))$ (counted with multiplicities)
coincide with the eigenvalues of 
$\alpha C_L+(1-\alpha)C_M$ (also counted with multiplicities), 
for every $\alpha \in 
\BR$. (This is a standard fact in 
the theory of matrix polynomials, see, for example,  
\cite{glr}.) 
Thus, $\alpha C_L+(1-\alpha)C_M$ has $n \ell$ real eigenvalues
for every $\alpha\in \BR$. 
Consider 
the homogeneous polynomial
of three real variables 
$\alpha,\beta, \gamma$:
$$ P(\alpha, \beta, \gamma):={\rm det}\, (\alpha C_L+\beta C_M- \gamma 
I_{n\ell}). $$
If $\alpha + \beta\neq 0$, the polynomial $P(\alpha, \beta, \gamma)$
(as a polynomial of $\gamma$) has $n\ell$ real roots (counted with 
multiplicities). If $\alpha+\beta=0$, then 
$$P(\alpha,\beta,\gamma)=\pm \gamma^{n(\ell-1)}\cdot 
{\rm det}\, (\alpha L_{\ell-1}+\beta M_{\ell-1}-\gamma I_n)$$
also has $n\ell$ real roots, by hypothesis of the theorem. 
Thus, $P(\alpha,\beta,\gamma)$ is hyperbolic in the direction 
of $(0,0,1)$, in the sense of the Lax conjecture, see \cite{L,
LPR}. 
By the main result of \cite{LPR} 
(the proof in \cite{LPR} is based on \cite{HV,V}), we have
$$ P(\alpha,\beta,\gamma)={\rm det}\, (\alpha A + \beta B - \gamma I)$$
for some real symmetric matrices $A$ and $B$.  
The direct statement of the theorem follows. 

To prove the converse statement, simply reverse the argument, taking into 
account that real symmetric matrices have all eigenvalues real. 
\eproof
\bigskip

For further development of the theory of hyperbolic 
polynomials of several variables and many applications, in particular, mixed 
determinants, see \cite{G}.

\section{Corollaries and applications}
\setcounter{equation}{0}

We start by recalling Obreschkoff's theorem
(see \cite{O,De}), 
which will be needed in the proof of the next corollary:
Two real scalar polynomials $f(z)$ and $g(z)$ of degrees 
$\ell$ and $\ell-1$, respectively, have the property that  
$f(z)+tg(z)$ has $\ell$ real roots (counted with multiplicities) for every 
real $t$ if and only if $f(z)$ and $g(z)$ have $\ell$
and $\ell-1$ real roots, respectively, and the roots of $f(z)$ and of $g(z)$
interlace (the cases when $f(z)$ or $g(z)$ have multiple roots and/or 
when $f(z)$ and $g(z)$ have common roots are not excluded here).

A proof of Obreschkoff's theorem can be given using the 
approach of Theorem \ref{th1}, as follows (below, we 
formulate 
Obreschkoff's theorem in a slightly different but equivalent form):
\begin{proposition} The following statements are
equivalent for  
scalar distinct monic relatively prime polynomials $f(z)$ and $h(z)$ of degree $\ell$:
\begin{itemize}
\item[(1)]  The polynomials $\alpha f(z)+ \beta h(z)$, $\alpha,\beta\in 
\BR$, $\alpha^{2} + \beta^{2} \neq 0$ ,
have all their roots real;
\item[(2)] The polynomials 
$\alpha f(z)+ (1-\alpha)h(z)$, $\alpha\in \BR$, have all their roots 
real;
\item[(3)] $f(z)$ has all its roots real, and the quotient 
$h(z)/f(z)$ has the form 
\begin{equation}\label{600}
 \frac{h(z)}{f(z)}=1+ \sum_{j=1}^p \frac{c_j}{z-\lambda_j}, 
\end{equation}
where $\lambda_j\in \BR, $ and the real numbers $c_j$ are all of the 
same sign;
\item[(4)]
Both $f(z)$ and $h(z)$ have $\ell$ distinct real roots and the 
roots of $f(z)$ and of $h(z)$ interlace. 
\end{itemize}\label{P:20} 
\end{proposition}

\proof 
(1) clearly implies (2) . (3) and (4) are equivalent : Indeed, if (3) 
or (4) holds true, then $f(z)$ has necessarily simple roots, 
and denoting the roots of $f(z)$ by 
$\lambda_1<\cdots < \lambda_p$, we see that 
in the representation (\ref{600}),
$$ {\rm sign }\, (c_j) ={\rm sign}\, \frac{h(z)}{f(z)} \qquad 
\mbox{as $z\rightarrow \lambda_j$, $z>\lambda_j$}, $$
whereas 
$$ -{\rm sign}\, (c_{j+1})= 
 {\rm sign}\, \frac{h(z)}{f(z)} \qquad 
\mbox{as $z\rightarrow \lambda_{j+1}$, $z<\lambda_{j+1}$}. $$
These equalities imply that that the roots of $f(z)$ and $h(z)$ 
interlace if and only if all the $c_j$'s are of the same sign. 

We next prove that (2) implies (3). Arguing as in the proof of Theorem 
\ref{th1}, we obtain that the characteristic polynomial
of $\alpha C_f+\beta C_h$ coincides with the 
characteristic polynomial of $\alpha A + \beta B$, for all 
$\alpha,\beta\in \BR$, 
where $A$ and $B$ are fixed (independent of $\alpha$ and $\beta$) distinct
real symmetric matrices.
Taking $\alpha - \beta=0$ we see that ${\rm rank}\, (A-B)\leq 1$. 
Since polynomials $f$ and $h$ are distinct 
it follow that $A=B\pm xx^T$ for some nonzero vector $x$. Now
$$ f(z)={\rm det}\, (zI-C_f)={\rm det}\, (zI- A),$$
$$ h(z)={\rm det}\, (zI- C_h)={\rm det}\, (zI-B),$$
and 
$$ \frac{h(z)}{f(z)}={\rm det}\, \left((zI-B)(zI-A)^{-1}\right) 
={\rm det}\, \left(I \pm xx^T(zI-A)^{-1}\right) =
1\pm x^T(zI-A)^{-1}x.$$ 
This reduces, upon applying a diagonalizing real orthogonal 
transformation $A\mapsto U^TAU$, and replacing $x$ with $U^Tx$,
to (\ref{600}) with 
the real numbers $c_j$ of the 
same sign, as required. 

Finally , let us prove prove the implication (3) $\ \Longrightarrow \ $ (1) .
This means to prove that if (3) holds then for any real $\gamma$
the equation $\frac{h(z)}{f(z)} + \gamma = 0$ does not have roots
with nonzero imaginary part . \\
Concider a complex number $z = a + bi$ , its real part $Re(z) = a$ ,
its imaginery part $Im(z) = b$ . If (3) holds then
$$
Im(\frac{h(z)}{f(z)}) = Im( \sum_{j=1}^p \frac{c_j}{z-\lambda_j} ) ,
$$
where $\lambda_j$ are real and the real numbers $c_j$ are all of the same sign .
Assume wlog that all $c_j$ are positive . As 
$$
Im((z-\lambda_j)^{-1}) = \frac{-b}{(a-\lambda_j)^{2} + b^{2} } ,
$$
thus we get that $Im(\frac{h(z)}{f(z)}) = -b \sum_{j=1}^p \frac{c_j}{(a-\lambda_j)^{2} + b^{2}}$ .\\
Therefore $Im(\frac{h(z)}{f(z)}) \neq 0$ if $Im(z) \neq 0$ . This means that
the equation $\frac{h(z)}{f(z)} + \gamma = 0$ does not have roots
with nonzero imaginary part for all real $\gamma$. \\
\eproof
\bigskip
   
We observe that checking condition (3) can be conveniently done using 
semidefinite programming. 
Indeed, let $h(z)$ and $f(z)$ be monic scalar polynomials with $f(z)$ 
having all roots real, and consider a minimal realization 
$$ \frac{h(z)}{f(z)}=1+ 
\widetilde{C}(zI-\widetilde{A})^{-1}\widetilde{B}, $$
where $\widetilde{C}$, $\widetilde{A}$, and $\widetilde{B}$ are real 
matrices. 
It is easy to see, using the uniqueness of a minimal realization up to a 
state isomorphism (similarity), that (3) holds, with the $c_j$'s 
positive, if and only if 
there exists a positive definite matrix $P$ such that 
$$ \widetilde{A}P=P\widetilde{A}^T, \qquad 
P\widetilde{C}^T=\widetilde{B}. $$
The latter problem is a semidefinite programming problem. \\
Another equivalent  semidefinite programming problem is based on the following
nice reformulation of Proposition 2.1 : \\
{ \it The conditions of Proposition 2.1 are equivalent to the existence of nonsingular real matrix $D$ such
that both $DC_{f}D^{-1}$ and $DC_{h}D^{-1}$ are real symmetric . } \\
(Our proof of this statement is essentially the same rank one pertubration argument as in
the proof of Proposition 2.1 ) . \\

This gives the following semidefinite programming problem : \\
is there exists a real positive definite $P \succ 0$ such that 
$$
PC_{f} = C_{f}^{T}P \mbox{\ and } PC_{h} = C_{h}^{T}P
$$

Our next corollary involves hyperbolic matrix polynomials.  
\begin{corollary}
Let $L(z)$ be a hyperbolic matrix polynomial. Then there exist 
$n\ell \times n\ell$ real symmetric matrices $A$ and $B$ such that 
the roots of ${\rm det}\, (L(z)+tL'(z))$ coincide with the eigenvalues of 
$A+tB$ $($multiplicities counted$)$, for every real number $t$.
Here, $L'(z)$ is the derivative of $L(z)$ with respect to $z$.
\end{corollary}

\proof.   
By Obreschkoff's theorem, the matrix polynomial $L(z)+tL'(z)$
 is hyperbolic for every real $t$. Now apply Theorem \ref{th1} with 
$M(z)=L(z)+L'(z).$
\eproof
\bigskip

Note that the condition (\ref{*}) implies (but is not equivalent to)
the condition that every convex combination of $L(z)$ and $M(z)$ is weakly 
hyperbolic. It turns out that the 
latter condition can be conveniently expressed
for hyperbolic matrix polynomials, which we will do next.

Let $L(\lambda)$ be a hyperbolic $n \times n$ matrix polynomial. 
For every $x\in \BC^n$, $\|x\|=1$, 
let 
$$ \lambda_1(x) \leq \lambda_2(x)\leq \cdots, \leq \lambda_{\ell}(x)$$
be the roots of equation (\ref{12}) arranged in the nondecreasing order.   
The sets 
$$ \Delta_j(L):=\{\, \lambda_j(x)\, | \, x\in \BC^n, \ \|x\|=1\},$$
called the {\em spectral zones} of $L(\lambda)$,  
are obviously closed intervals on the real line:
$$  \Delta_j(L)=[\delta^-_j(L), \delta^+_j(L)], \qquad j=1,2, \ldots, \ell.$$
A basic result in  
the theory of hyperbolic matrix and operator polynomials 
(\cite[Theorem 31.5]{3}, for example), states that 
two spectral zones either are disjoint, or have only one 
point in common.
 
\begin{proposition} Let $L(\lambda)$ and
$M(\lambda)$ be two hyperbolic matrix polynomials of degree $\ell$.
Then every convex combination $
 \alpha L(z)+(1 -\alpha)M(z)$, $0\leq \alpha\leq 1$, is 
hyperbolic if and only if their spectral zones satisfy the inequalities
$$ \max\{\delta^+_j(L),\delta^+_j(M)\}
\leq \min\{\delta^-_{j+1}(L),\delta^-_{j+1}(M)\}, \qquad j=1,\ldots, 
\ell-1.$$
\end{proposition}

For the proof apply \cite[Theorem 2.1]{De}; this theorem gives 
necessary and 
sufficient conditions for all linear combinations of two 
given scalar polynomials to be hyperbolic.  

Using Theorem \ref{th1} and inequalities for eigenvalues of 
real symmetric matrices (see for example \cite{LO}), one can derive 
inequalities for eigenvalues of weakly hyperbolic matrix polynomials.
We illustrate this for
the case of the Horn inequalities. 
For a Hermitian $m \times m$ matrix $X$, we write its eigenvalues (repeated 
according to 
their multiplicities) in a non-decreasing order:
$$ \lambda_1(X)\leq\lambda_2(X)\leq \cdots \leq \lambda_m(X). $$   
An ordered triple $(U,S,T)$ of nonempty subsets of $\{1,2,\ldots, m\}$ 
is said to be a  {\em Horn triple} (with respect to $m$) if the 
cardinalities of  
$U$, $S$, and $T$ are the same, and the 
{\em Horn inequalities} 
$$ \sum_{i\in U} \lambda_i(X+Y)\leq
\sum_{j\in S} \lambda_j(X) + \sum_{k\in T} \lambda_k(Y) $$
hold true for every pair of Hermitian $m \times m$ matrices $X$ and $Y$.
A description of all Horn triples is known \cite{Kl, Kl1}; see 
also 
the surveys \cite{Fu, Bh}. 
For a weakly hyperbolic $n \times n$ matrix polynomial $L(z)$ of degree 
$\ell$,
we arrange the roots of ${\rm det}\, (L(z))$ in the non-decreasing order:
$$ d_1(L)\leq d_2(L) \leq \ldots \leq d_{n\ell}(L). $$
Let $T = \{1 \leq i_1 < i_2 < ...< i_m \leq n\ell \} \subset \{1,2,...,n\ell \}$ , define
$\bar{T} = \{n\ell -i_m < n\ell -i_{m-1} < ...<  n\ell -i_1 \} $ .

\begin{theorem}
Let $L(z)$ and $M(z)$ be monic $n\times n$ matrix polynomials satisfying the 
hypotheses of Theorem $\ref{th1}$. Then for every Horn triple $(U,S,T)$ 
with respect to $n\ell$, and for every $\alpha\in \BR$, the inequality
$$
 {\displaystyle \sum_{i\in U} d_i\left(\alpha L+(1-\alpha)M\right) \leq
\alpha
\left(\sum_{j\in S_{\alpha} } d_j(L) \right)+ (1-\alpha)\left(\sum_{k\in T_{1-\alpha}} d_k(M)\right)} 
$$
holds true. Here $S_{\alpha} = S$ if $\alpha \geq 0$ and $S_{\alpha} = \bar{S}$ if $\alpha < 0$.
\end{theorem}

\proof Let $A$ and $B$ be as in Theorem \ref{th1}.  
Then we have, using Theorem \ref{th1} and the Horn inequalities:
\begin{eqnarray*}
 {\displaystyle \sum_{i\in U} d_i\left(\alpha L+(1-\alpha)M)\right)
} & = & {\displaystyle \sum_{i\in U} \lambda_i\left(\alpha A+(1-\alpha )B\right)
}\\ &\leq &
{\displaystyle \sum_{j\in S} \lambda_j\left(\alpha A\right) + \sum_{k\in T} 
\lambda_k\left((1-\alpha )B\right)} \\ 
&=& {\displaystyle \alpha \left(\sum_{j\in S_{\alpha}} \lambda_j(A)\right) + 
(1-\alpha)\left(\sum_{k\in 
T_{1-\alpha}} \lambda_k(B)
\right)} \\ &=& {\displaystyle  
\alpha\left(
\sum_{j\in S_{\alpha}} d_j(L)\right) + (1-\alpha)\left(\sum_{k\in T_{1-\alpha}} d_k(M)\right)},
\end{eqnarray*}
and the proof is complete.

\end{document}